\documentclass[11pt]{article}
\newcommand{\E}{\bf{E}}
\newcommand{\D}{\bf{D}}
\newcommand{\N}{\bf{N}}
\newcommand{\Z}{\bf{Z}}
\newcommand{\R}{\bf{R}}
\newtheorem{conjecture}{Conjecture}

\newtheorem{corollary}[conjecture]{Corollary}

\newtheorem{lemma}[conjecture]{Lemma}

\newtheorem{theorem}[conjecture]{Theorem}

\begin{document}

\newenvironment{definition}{\noindent {\bf Definition.} \rm}{\hfill \\ \vspace{.15cm}}

\newenvironment{remark}{\noindent {\bf Remark.}  \rm}{\hfill \\ \vspace{.15cm}}
\newenvironment{proof}{\noindent {\bf Proof} \hspace{.1cm}}{\hfill $\diamondsuit$ \\ \vspace{.15cm}}

\title{The Weil-Petersson geometry of the five-times punctured sphere}
\author{Javier Aramayona}
\date{May 2004 (Revised, February 2005)}

\maketitle

\begin{abstract}
We give a new proof that the completion of the Weil-Petersson
metric on Teichm\"uller space is Gromov-hyperbolic if the surface
is a five-times punctured sphere or a twice-punctured torus. Our
methods make use of the synthetic geometry of the Weil-Petersson
metric.
\end{abstract}

\section{Introduction}

The large scale geometry of Teichm\"uller space has been a very
important tool in different aspects of the theory of hyperbolic
$3$-manifolds. Within this context, a natural question to ask is
whether Teichm\"uller space, with a given metric, is hyperbolic in
the sense of Gromov (or Gromov hyperbolic, for short). In general,
the answer is negative. In their paper \cite{BrFa}, the authors
prove the following: if $\Sigma$ is a surface of genus $g$ and
with $p$ punctures, with $3g-3+p> 2$, then the Teichm\"uller space
of $\Sigma$, endowed with the {\it{Weil-Petersson metric}}, is not
Gromov hyperbolic. In the case when $3g-3+p = 2$ (that is, when
the surface is a sphere with five punctures or a twice-punctured
torus) they show that the Weil-Petersson Teichm\"uller space is
Gromov hyperbolic. The proof makes reference to very deep results
by Masur and Minsky on the Gromov hyperbolicity of the {\it curve
complex}. The aim of this paper if to give a direct proof of the
Gromov hyperbolicity of the Weil-Petersson Teichm\"uller space in
the case $3g-3+p = 2$.

The objective of this article is to show the following result:

\begin{theorem} \label{main result}
If $\Sigma$ is the five-times punctured sphere or the twice
punctured torus, then $\overline {\cal{T}}_{WP}(\Sigma)$ is Gromov
hyperbolic.
\end{theorem}

Brock \cite{Br} showed that, for every hyperbolic surface $S$,
${\overline {\cal{T}}}_{WP}(\Sigma)$ is quasiisometric to the
pants complex (see \cite{Br} for definitions)
${\cal{C}}_P(\Sigma)$ of the surface $\Sigma$. Since Gromov
hyperbolicity is a quasiisometry invariant, we get the following
result:

\begin{corollary}
If $\Sigma$ is the five-times punctured sphere or the twice
punctured torus, then the pants complex ${\cal{C}}_P(\Sigma)$ of
$\Sigma$ is Gromov hyperbolic.
\end{corollary}

\begin{remark}
Behrstock \cite{Be} has given a direct proof of this last result,
using the combinatorial structure of the pants complex.
\end{remark}

Let $\Sigma$ be a hyperbolic surface of genus $g$ and with $p$
punctures and let ${\cal{T}}(\Sigma)$ be the Teichm\"uller space
of $\Sigma$. The Weil-Petersson metric on ${\cal{T}}(\Sigma)$ is a
non-complete metric of negative sectional curvature.

\begin{remark}
We note that Theorem \ref{main result} cannot be obtained as a
consequence of the negative sectional curvatures of the
Weil-Petersson metric, since these curvatures have been shown (see
\cite{Hu}) not to be bounded away from zero.
\end{remark}

Let ${\cal{T}}_{WP}(\Sigma)$ denote the Teichm\"uller space of
$\Sigma$ endowed with the Weil-Petersson metric. The metric
completion ${\overline {\cal{T}}}_{WP}(\Sigma)$ of
${\cal{T}}_{WP}(\Sigma)$ is the {\it{augmented Teichm\"uller space
}} (see \cite{Ma}), i.e. the set of marked metric structures on
$\Sigma$ with nodes on a (possibly empty) collection of different
homotopy classes of essential simple closed curves on $\Sigma$ .
Here, a curve is {\it{essential}} if it is not null-homotopic nor
homotopic to a puncture. From now on we will refer to a homotopy
class of essential simple closed curves on $\Sigma$ simply as a
{\it {curve}}, unless otherwise stated.

Let ${\cal{C}} = {\cal{C}}(\Sigma)$ be the {\it{curve complex}} of
$\Sigma$, as defined in \cite{Har}. Recall that this is a
finite-dimensional simplicial complex whose vertices correspond to
(homotopy classes of non-trivial, non-peripheral) curves on
$\Sigma$, and that a subset $A \subseteq V({\cal{C}})$ spans a
simplex in $\cal{C}$ if the elements of $A$ can be realised
disjointly on $\Sigma$. Note that inclusion determines a partial
order on the set of simplices of ${\cal{C}}$. Following \cite{Wo},
we can define a map $\Lambda: {\overline {\cal{T}}}_{WP}(\Sigma)
\rightarrow {\cal{C}} \cup \{\emptyset\}$ that assigns, to a point
in $u \in {\overline{\cal{T}}}_{WP}(\Sigma)$, the (possibly empty)
collection of different curves on $\Sigma$ on which $u$ has nodes.
The space ${\overline {\cal{T}}}_{WP}(\Sigma)$ is the union of the
level sets of $\Lambda$. Then, ${\overline
{\cal{T}}}_{WP}(\Sigma)$ has the structure of a {\it{stratified}}
space, where the level sets of $\Lambda$ are the {\it{strata}}
(observe that $\Lambda^{-1}(\{\emptyset\}) =
{\cal{T}}_{WP}(\Sigma)$ and that two strata intersect over a
stratum if at all). We will refer to $\Lambda^{-1}(\Lambda(u))$ as
the {\it{stratum containing $u$}} and we will say that it has
label $\Lambda(u)$. The strata of ${\overline
{\cal{T}}}_{WP}(\Sigma)$ are isometric embeddings of products of
lower dimensional Teichm\"uller spaces (which come from
subsurfaces of $\Sigma$) with their corresponding Weil-Petersson
metric. It is clear that the stratum with label a collection of
curves that determine a pants decomposition on $\Sigma$ consists
of only one point in ${\overline {\cal{T}}}_{WP}(\Sigma)$. Let
${\rm{Mod}}(\Sigma)$ be the mapping class group of $\Sigma$, i.e.
the group of self-homeomorphisms of $\Sigma$ up to homotopy. It is
known (see \cite{Ab}) that ${\rm{Mod}}(\Sigma)$ acts cocompactly
on ${\overline {\cal{T}}}_{WP}(\Sigma)$. The space ${\overline
{\cal{T}}}_{WP}(\Sigma)$ is not locally compact: a point in
${\overline {\cal{T}}}_{WP}(\Sigma) \setminus
{\cal{T}}_{WP}(\Sigma)$ does not admit a relatively compact
neighbourhood. Indeed, let $u$ be a point in ${\overline
{\cal{T}}}_{WP}(\Sigma) \setminus {\cal{T}}_{WP}(\Sigma)$, which
corresponds to a surface with a nodes on the simple closed curve
$\alpha$ (and possibly more), and consider the Dehn twist
$T_{\alpha}$ along $\alpha$. Then the $T_{\alpha}$ orbit of any
point lies in every neighbourhood of $u$ (see \cite{Wo}).
Nevertheless, individual strata are locally compact, since they
are (products of) lower dimensional Teichm\"uller spaces.

The following result summarises some deep and remarkable facts
about the geometry of the Weil-Petersson metric on Teichm\"uller
space (see below for the relevant definitions). Part $(i)$ is due
to S. Yamada \cite{Ya}; $(ii)$ is due to Wolpert \cite{Wo} and
$(iii)$ is due to Daskalopoulos and Wentworth \cite{DW}. Let us
note that all these results rely on previous work by Wolpert on
the Weil-Petersson metric.

\clearpage

\begin{theorem}[\cite{DW}, \cite{Wo}, \cite{Ya}]
Let $\Sigma$ be a surface of hyperbolic type and let
${\overline{\cal{T}}}_{WP}(\Sigma)$ be the completion of the
Weil-Petersson metric on ${\cal{T}}_{WP}(\Sigma)$. Then,

\begin{enumerate}
\item The space ${\overline {\cal{T}}}_{WP}(\Sigma)$ is a
${\rm{CAT}}(0)$ space.

\item The closure of a stratum in ${\overline
{\cal{T}}}_{WP}(\Sigma)$ is convex and complete in the induced
metric.

\item The open geodesic segment $[u,v] \setminus \{u,v\}$ from $u$
to $v$ lies in the stratum with label $\Lambda(u) \cap
\Lambda(v)$.
\end{enumerate}
\end{theorem}

Let $X_0$ be the Teichm\"uller space of the five-times punctured
sphere $\Sigma_{0,5}$ endowed with the Weil-Petersson metric and
let $X$ be its metric completion. We will write $X_F = X \setminus
X_0$. Since a pants decomposition of $\Sigma_{0,5}$ corresponds to
two disjoint curves on $\Sigma_{0,5}$, we get that a stratum in
$X$ has label $\alpha$ or $\alpha\beta$, where $\alpha$ and
$\beta$ are disjoint curves on $\Sigma_{0,5}$. The closure of a
stratum of the type $S_{\alpha}$ is given by
$\overline{S}_{\alpha} = S_{\alpha} \cup (\bigcup_{\beta \in
{\cal{B}}} S_{\alpha\beta})$, where $\cal B$ is the set of curves
that are disjoint from $\alpha$. We observe that any curve
$\alpha$ separates $\Sigma_{0,5}$ into two subsurfaces, namely a
four-times punctured sphere and a three-times punctured sphere.
Then $S_{\alpha}$ corresponds to the Teichm\"uller space of the
four-times punctured sphere, since the Teichm\"uller space of the
other subsurface is trivial. Also, recall that a stratum of the
form $S_{\alpha\beta}$ is a single point in $X$.

We prove Theorem \ref{main result} in the case where the surface
is a sphere with five punctures. Using the techniques from
Sections $2$ and $3$ we immediately get the result for the
twice-punctured torus. The only difference between the two cases
is the nature of the strata in $\overline{T}_{WP}(\Sigma)$ which
arise from pinching a single curve on the surface. In the first
case, these strata correspond to the Teichm\"uller space of a
four-times punctured sphere; in the second, they correspond either
to the Teichm\"uller space of a four-times punctured sphere or to
the Teichm\"uller space of a once-punctured torus, depending on
whether the curve giving rise to such a stratum separates the
surface. In both cases, $\overline{T}_{WP}(\Sigma)$ has the same
structure (as a stratified space).

\noindent {\it {Acknowledgements.}} I would like to thank Brian
Bowditch for introducing me to the problem and for the interesting
conversations we had about the topic. I would also like to thank
the referee for his/her useful comments.

\section{Preliminaries}

In order to give a proof of Theorem \ref{main result} we will have
to make use of some geometric properties of ${\rm{CAT}}(0)$
spaces.

\noindent Let us begin by recalling the definition of a
${\rm{CAT}}(0)$ space. Let $Y$ be a {\emph{geodesic metric
space}}, that is, a space in which every two points in the space
can be connected by a path which realises their distance; such a
path is called a {\emph{geodesic}} between the two points. By a
{\emph{triangle}} in $Y$ we will mean three points $x_1,x_2,x_3
\in Y$, the {\emph{vertices}} of $T$, and three geodesics
connecting them pairwise. We will write $[x_i,x_j]$ for the
geodesic side of $T$ with endpoints $x_i$ and $x_j$. Throughout
this article we will denote the euclidean plane by $\E^2$ and the
euclidean distance in $\E^2$ by $d_e$.

\begin{definition}
Let $Y$ be a geodesic metric space and let $T$ be  triangle in $Y$
with vertices $x_1,x_2,x_3$. A \em{comparison triangle} for $T$ in
$\E^2$ is a geodesic triangle in $E^2$ with vertices
$\overline{x}_1$, $\overline{x}_2$ and $\overline{x}_3$ such that
$d(x_i,x_j)=d_e(\overline{x}_i, \overline{x}_j)$ for all
$i,j=1,2,3$. Given a point $p \in [x_i,x_j]$, for some $i,j=1,2,3$
distinct, a \em{comparison point} $\overline{p}$ for $p$is a point
$\overline{p} \in [\overline{x}_i, \overline{x}_j]$ such that
$d(x_i,p)=d_e(\overline{x}_i,\overline{p}$ and
$d(x_j,p)=d_e(\overline{x}_j,\overline{p}$.
\end{definition}

\begin{definition}
We say that the triangle $T$ satisfies the \em{${\rm{CAT}}(0)$
inequality} if for any points $p \in [x_i,x_j]$ and $q \in
[x_j,x_k]$, for $i,j,k  = 1, 2, 3$ distinct, we have that $d(p,q)
\leq d_e( \overline{p}, \overline{q})$, where $\overline{p}$ and
$\overline{q}$ are comparison points for $p$ and $q$,
respectively, in the comparison triangle $\overline{T}$ for $T$.
We will say that the space $X$ is a \em{${\rm{CAT}}(0)$ space} if
every triangle in $X$ satisfies a ${\rm{CAT}}(0)$ inequality.
\end{definition}

The next three results about ${\rm{CAT}}(0)$ spaces are
well-known; they will be crucial in our main argument. For a proof
see, for instance \cite{BriHa}.

\begin{theorem} \label{convexity}
If $Y$ is a ${\rm{CAT}}(0)$ space, then the distance function on
$Y$ is convex along geodesics, that is, if $\sigma,\sigma':[0,1]
\rightarrow Y$ are geodesics in $Y$ parametrised proportional to
arc-length, then

$$d(\sigma(t),\sigma'(t))\leq t d(\sigma(0),\sigma'(0)) + (1-t) d(\sigma(1),\sigma'(1)),$$

\noindent for all $t \in [0,1]$.
\end{theorem}

\begin{corollary}
Every ${\rm{CAT}}(0)$ space is uniquely geodesic.
\end{corollary}

\begin{theorem} \label{projection}
Let $Y$ be a ${\rm{CAT}}(0)$ space and let $C$ be a complete
convex subset of $Y$. Given $x\in Y$ there exists a point $\pi(x)
\in C$ such that $d(x,\pi(x))={\rm{inf}}_{c\in C} d(x,c)$.
Moreover, the map $\pi:Y \rightarrow C$ is distance
non-increasing, that is, for all $x,y \in Y$ we have that $d(x,y)
\leq d_C(\pi(x), \pi(y))$, where $d_C$ denotes the subspace
metric.
\end{theorem}

When interested in the large-scale geometry of a ${\rm{CAT}}(0)$
space, a natural question to ask is what are the obstructions for
such a space to be Gromov hyperbolic. An answer to this question
was given by Bowditch \cite{Bo} and Bridson \cite{Bri} separately,
generalising a result announced by Gromov. They showed the
following.

\begin{theorem} [\cite{Bo}, \cite{Bri}]
Let $Y$ be complete, locally compact ${\rm{CAT}}(0)$ space which
admits a cocompact isometric group action. Then, either $Y$ is
Gromov hyperbolic or else it contains a totally geodesic embedding
of a eulidean plane.
\end{theorem}

Let $X$ be the completion of the Weil-Petersson metric on the
Teichm\"uller space of the five-times punctured sphere. We are
going to show, using some of the techniques in \cite{Bo}, that if
$X$ is not Gromov hyperbolic then there exists an isometrically
embedded euclidean disc in one of the strata of $X$, which is
impossible since all the sectional curvatures of $X$ are strictly
negative. More specifically, we are going to construct, for each
$n\in\N$, a map $\phi_n:\D \rightarrow X$, where $\D$ is the unit
disc in $\E^2$, such that

$$ \lambda_n d_e(x,y) \leq d(\phi_n(x),\phi_n(y))\leq d_e(x,y),$$

\noindent for all $x,y \in \D$ , and where $\lambda_n \in (0,1)$
with $\lambda_n \rightarrow 1$ as $n \rightarrow \infty$. If the
space $X$ were locally compact, one could take the limit of
$phi_n$ as $n$ tends to infinity, obtaining in this way an
isometric embedding of the euclidean unit disc in $X$. Even though
$X$ is not locally compact, we are able to use the structure of
$X$ as a stratified space and the fact that individual strata are
locally compact to obtain such an isometrically embedded disc. We
note that we will not use the cocompact isometric action of the
mapping class group on $X$ to construct these maps; this action
will only be used in the arguments in the next section. \noindent
The construction of the maps $\phi_n$ is totally analogous to the
one in \cite{Bo}, where Bowditch shows the following result (we
remark that the results in \cite{Bo} are more general than the
ones we present here).

\begin{lemma}[\cite{Bo}]
Let $Y$ be a ${\rm{CAT}}(0)$ space. Given $n \in \N$ and $\epsilon
>0$, there exists a
map $\phi_{n,\epsilon}: ([-n,n] \cap \Z)^2 \rightarrow Y$ such
that
$$ \mid i - i' \mid - \epsilon \leq d(\phi_{n,\epsilon}(i,j),\phi_{n,\epsilon}(i',j)) \leq \mid
i - i' \mid \eqno(1)$$ and $$ \mid j - j' \mid - \epsilon \leq
d(\phi_{n,\epsilon}(i,j),\phi_{n,\epsilon}(i,j')) \leq \mid j - j'
\mid \eqno(2)$$
\end{lemma}

\begin{remark}
For the sake of completeness, we now give a brief account on
Bowditch's construction. Let $n\in\N$ and $\epsilon>0$ and let $q$
ne a natural number bigger than $2n^2 /over \epsilon$. Let
$\sigma:[0,qn] \rightarrow Y$ be a geodesic segment in $Y$ and let
$y$ be a point in $Y$ at distance at least $n$ from $\alpha$. Let
$\tau_i:[0,d(y,\sigma(i)] \rightarrow Y$ be the uniue geodesic
from $y$ to $\sigma(i)$, for all $i \in [0,qn] \cap \N$. Bowditch
then sets, for all $i,j=0, dots, n$ and for all $p=0, \dots, q-1$,
$phi_p(i,j)=\tau_{pn+i}(j)$ and shows that there exists a number
$p=0, \dots n-1$ satisfying (1) and (2).
\end{remark}

Bowditch then shows that if the space $Y$ is, in addition, not
Gromov hyperbolic then the maps $\phi_n$ satisfy the following
non-degeneracy condition.

\begin{lemma}
Let $Y$ be a ${\rm{CAT}}(0)$ space and suppose $Y$ is not Gromov
hyperbolic. Then the maps $\phi_n$ in the result above satisfy, in
addition, that
$$d(\phi_{n,\epsilon}(i,j),\phi_{n,\epsilon}(i+1,j+1)) \geq 1/2 \eqno(3)$$ and
$$d(\phi_{n,\epsilon}(i,j+1), \phi_{n,\epsilon}(i+1,j)) \geq 1/2 \eqno(4)$$ for all $n \in \N$.
\end{lemma}

We now give an extension, using standard arguments about
${\rm{CAT}}(0)$ spaces, to Bowditch's construction. It will play a
central role in the proof of Theorem \ref{main result}.

\begin{lemma}
Let $Y$ be a ${\rm{CAT}}(0)$ space and suppose that $Y$ is not
Gromov hyperbolic. Fix a number  $N \in \N$ and consider $K=[-N,N]
\times  [-N, N] \subseteq \R^2$, endowed with the euclidean
metric. Then there exists a sequence $(\phi_n)_{n\in\N}$ of maps
$\phi_n:K \rightarrow Y$ such that

$$\lambda_n d_e(x,y) \leq d(\phi_n(x), \phi_n(y)) \leq d_e(x,y),$$

\noindent where $\lambda_n \in (0,1)$ for all $n\in\N$ and
$\lambda_n \rightarrow 1$ as $n \rightarrow \infty$.
\end{lemma}

\begin{proof}
Let $(\epsilon_n)_{n\in\N}$ be a sequence of positive real numbers
tending to zero and let $\phi_n = \phi_{n,\epsilon_n}$ be the map
described above. First, we are going to extend the map $\phi_n$ to
$K$ as follows:

\noindent Let $\sigma_n: [0,1] \rightarrow Y$  be the unique
geodesic segment in $Y$ from $\phi_n(0,0)$ to $\phi_n(1,0)$.
Similarly, let $\sigma'_n: [0,1] \rightarrow Y$ be the geodesic
segment in $Y$ from $\phi_n(0,1)$ to $\phi_n(1,1)$. Here, and from
now on, we assume that all the geodesics are parametrised
proportional to arc-length. For $t \in [0,1]$, let $\tau^n_t:[0,1]
\rightarrow Y$ be the unique geodesic in $Y$ connecting
$\sigma(t)$ and $\sigma'(t)$. We set $\phi_n(t,s):= \tau^n_t(s)$,
for all $t,s \in [0,1]$. We can extend $\phi_n$ to $K$ by
performing this construction on each square $[i,i+1] \times
[j,j+1]$ for all $i,j \in [-N,N] \cap \Z$, provided $n$ is large
enough so that $\phi_n$ is defined on the shole of $K$. Note that,
up to extracting a subsequence, we can assume that this is the
case.

\smallskip

Since the distance function on $Y$ (Theorem \ref{convexity}) is
convex along geodesics we get, from inequality (2), that

$$d(\phi_n(t,j),\phi_n(t,j+1)) \leq 1,$$

for all $t \in [-N,N]$, $j \in [-N,N] \cap \Z$ and $n\in\N$. The
convexity of the distance function on $Y$ yields that the real
function $[t \rightarrow d(\phi_n(t,j),\phi_n(t,j+1))]$ tends, as
$n$ grows, to a real function which is convex (and also bounded,
from (2)). But a real function which is convex and bounded must be
constant and therefore $d(\phi_n(t,j),\phi_n(t,j+1) \rightarrow 1$
as $n \rightarrow \infty$, from $(1)$.

\smallskip

Note that, from $(2)$ and the properties of the extension $\phi_n$
to $K$, we can deduce that the images of a vertical segment of the
form $\{i\} \times [-N,N]$ under $\phi_n$, where  $i \in [-N,N]
\cap \Z$, get {\emph{arbitrarily close}} to being geodesic in $Y$
as $m$ grows; by this we mean that the Hausdorff distance between
$\phi_n(\{i\} \times [-N,N])$ and the unique geodesic in $Y$ with
the same endpoints tends to 0 as $n$ tends to infinity. By a
totally analogous convexity argument we obtain that the images of
any two vertical lines in $K$ under $\phi_n$ get arbitrarily close
to being parallel geodesics as $n \rightarrow \infty$. Note that
this implies that, up to taking a subsequence, the maps $\phi_n$
are injective on $K$.

\smallskip

Let $0 \leq t,t',s,s' \leq 1$. Again from the convexity of the
distance function we obtain that

\begin{eqnarray}
d(\tau^n_t(s), \tau^n_{t'}(s)) &\leq& s d(\tau^n_t(0),
\tau^n_{t'}(0)) + (1-s) d(\tau^n_t(1), \tau^n_{t'}(1)) \nonumber
\\ & =& s d(\sigma_n(t), \sigma_n(t')) + (1-s) d(\sigma_n'(t),
\sigma_n'(t')) {} \nonumber \\  & =& \mid t -t' \mid {}. \nonumber
\end{eqnarray}

\noindent From this inequality and the fact that $\tau^n_t: [0,1]
\rightarrow Y$ is a geodesic paremetrised proprotional to arc
length for every $t \in [0,1]$, we deduce that

\begin{eqnarray}
d(\phi_n(t,s), \phi_n(t',s')) &=& d(\tau^n_t(s), \tau^n_{t'}(s'))
 \nonumber \\ & =& d(\tau^n_t(s), \tau^n_{t}(s')) +
d(\tau^n_t(s'), \tau^n_{t'}(s')) \nonumber \\ & \leq& \mid s - s'
\mid + \mid t - t' \mid {}, \nonumber
\end{eqnarray}

\noindent and thus $\phi_n$ is continuous on $K$ for all $n\in\N$.

\smallskip

Recall that the images of vertical lines in $K$ under the maps
$\phi_n$ get arbitrarily close to being parallel geodesis in $Y$
as $n$ grows. Consider now a horizontal line $[-N,N] \times
\{s_0\}$ in $K$  and let $\rho_n:[0,1] \rightarrow Y$ be the
geodesic in $Y$ between $\phi_n(0,s_0)$ and $\phi_n(1,s_0)$. We
know, from the discussion above, that $d(\sigma_n(t),
\sigma'_n(t)) \rightarrow 1$ as $n \rightarrow \infty$ for all $t
\in [0,1]$. From the convexity of the distance function we get
that $d(\sigma_n(t), \rho_n(t)) \rightarrow s_0$ and
$d(\sigma'_n(t), \rho_n(t)) \rightarrow 1 - s_0$ as $n \rightarrow
\infty$. But $d(\phi_n(t,s_0), \sigma_n(t)) \rightarrow s_0$ and
$d(\phi_n(t,s_0), \sigma'_n(t)) \rightarrow 1 - s_0$  as $n
\rightarrow \infty$, since $[s \rightarrow \phi_n(t,s)]$ is
geodesic for a fixed $t \in [0,1]$. Thus the images of a
horizontal line in $K$ under $\phi_n$ are paths that get
arbitrarily close to being geodesics $n \rightarrow \infty$. By a
similar argument we get that the same holds for the images of a
horizontal line in $K$ under $\phi_n$.

\smallskip

Let $d_n$ be the pull-back metric on $K$ determined by $\phi_n$.
Since $d_n((-N,-N),(N,N))$ is bounded above and below (note that,
in particular, it is bounded below away from 0, from (1)) we can
assume that $d_n((-N,-N),(N,N)) \rightarrow a> 0$, up to
extracting a subsequence. Also recall, from inequalities (1) and
(2), that $d_n((-N,-N),(-N,N)) \rightarrow 2N$ and
$d_n((-N,-N),(N,-N)) \rightarrow 2N$ as $n \rightarrow \infty$.

\noindent Let $f :\R^2 \rightarrow \R^2$ be an affine map such
that $f(t,0) = (t,0)$, for all $t \in \R$, and 
$d_e(f(-N,-N),f(N,N)) = a$. By considering $\phi_n \circ f^{-1}$,
which we denote again by $\phi_n$ abusing notation, we deduce that
the maps $\phi_n$, restricted to $K$, satisfy that

$$d_e(x,y) - \epsilon_n \leq d(\phi_n(x),\phi_n(y)) \leq d_e(x,y),$$

\noindent for all $x,y \in K$. Consider the restriction of
$\phi_n$ to the unit disc $\D$ in $\R^2$. It follows immediately
from the last inequality that there exists a sequence
$(\lambda_n)_{n\in\N}$ of real numbers, with $\lambda_n =
\lambda_n(\epsilon_n)$ and $\lambda_n \in (0,1)$, such that

$$\lambda_n d_e(x,y) \leq d(\phi_n(x),\phi_n(y)) \leq d_e(x,y),$$

\noindent as desired.
\end{proof}

\section {Proof of the Theorem}

Let us begin with some technical results that will be important to
prove Theorem \ref{main result}.

\begin{remark} Let $\Sigma$ be a hyperbolic surface. In \cite{Wo},
Corollary $21$, Wolpert shows that for a stratum $S$ defined by
the vanishing of the length sum $l = l_1 + \dots l_n$, where $l_i$
corresponds to the length of some curve $\alpha_i$ on the surface
for $i = 1, \dots, n$, the distance to the stratum is given
locally as $d(p, S)= (2\pi l)^{1/2} + O(l^2)$. In particular, if
$(u_m)_{m\in\N}$ is a sequence in $X$ such that $d(u_m,S)
\rightarrow 0$ as $m \rightarrow \infty$ we get that
$l_{u_m}(\alpha_i) \rightarrow 0$ as $m \rightarrow \infty$, for
all $i = 1, \dots, n$, where $l_{u_m}(\alpha)$ denotes the length
of $\alpha$ in $u_m$.
\end{remark}

\begin{lemma}
Let $S_{\alpha}$ and $S_{\beta}$ be two different strata in $X_F$
such that $\overline{S}_{\alpha} \cap \overline{S}_{\beta} \neq
\emptyset$. Let $(u_n)_{n\in\N}$ and $(v_n)_{n\in\N}$ be sequences
in $\overline{S}_{\alpha}$ and $\overline{S}_{\beta}$,
respectively, such that $d(u_n,v_n) \rightarrow 0$ as $n
\rightarrow \infty$. Then $d(u_n, S_{\alpha \beta}) \rightarrow 0$
and $d(v_n, S_{\alpha \beta}) \rightarrow 0$ as $n \rightarrow
\infty$.
\end{lemma}

\begin{proof}
Suppose, for contradiction, that the result is not true. Let
$\Gamma$ be the mapping class group of $\Sigma_{0,5}$. Since the
action of $\Gamma$ on $X$ is cocompact there is a compact subset
$Z$ of $X$ and a sequence $(\gamma_n)_{n\in\N}$ in $\Gamma$ such
that $\gamma_n(u_n)\in Z$ for all $n\in\N$. Note that $d(u_n,v_n)
\rightarrow 0$ as $n \rightarrow \infty$ if and only if
$d(\gamma_n(u_n),\gamma_n(v_n)) \rightarrow 0$ as $n \rightarrow
\infty$, since $\Gamma$ acts isometrically on $X$ and therefore
$d(\gamma_n(u_n),\gamma_n(v_n))=d(u_n,v_n)$ for all $n\in\N$. We
replace $u_n$ and $v_n$ by $\gamma_n(u_n)$ and $\gamma_n(v_n)$,
respectively; abusing notation we denote the new points by $u_n$
and $v_n$ again. Then, up to extracting a subsequence, we get that
the sequence $(u_n)_{n\in\N}$ converges to a point $u_0 \in X$,
since $Z$ is compact. Since $u_n \in \overline{S}_{\alpha}$, then
$l_{u_n}(\alpha) = 0$ for all $n\in\N$. Similarly, we have that
$l_{v_n}(\beta) = 0$ for all $n$. From the remark above and since
$d(u_n, S_{\alpha \beta}) >K$ for infinitely many $n$, it follows
that there must exist a constant $ L = L(K)>0$ such that
$l_{u_n}(\beta) > L$, for infinitely many $n\in\N$. Therefore, for
the simple closed curve $\beta$ we have that $l_{u_n}(\beta)>L>0$
and $l_{v_n}(\beta)= 0$, for infinitely many $n$, which is
impossible since we know from the hypotheses that $d(u_n,v_n)
\rightarrow 0$ as $n \rightarrow \infty$
\end{proof}

As a trivial consequence, we get the following corollary:

\begin{corollary}
Given $r>0$ and a stratum of the form $S_{\alpha\beta}$, let $B(r)
= B(S_{\alpha \beta}, r)$ be the ball of radius $r$ around
$S_{\alpha \beta}$. Then, there exists $D = D(r)>0$ such that
$d(u,v)>D$, for all $u \in \overline{S}_{\alpha} \setminus (B(r)
\cap \overline{S_\alpha})$ and $v \in \overline{S}_{\beta}
\setminus (B(r) \cap \overline{S_\beta})$.
\end{corollary}

\begin{proof}[of Theorem \ref{main result}]

Suppose, for contradiction, that $X$ is not Gromov hyperbolic.
Since $X$ is a ${\rm{CAT}}(0)$ space, we can construct a sequence
of maps $(\phi_n:\D \rightarrow X)_{n\in\N}$ as described in the
previous section. Using the structure of $X$ as a stratified
space, there is a stratum $S \subseteq X$ such that $\phi_n(\D)
\subseteq \overline{S}$ (note that $S$ may have label $\emptyset$,
so that $\overline{S} = X$. Also, $S$ cannot have a pants
decomposition as label, since recall that such a stratum consists
of only one point). The stratum $S$ corresponds to the
Weil-Petersson Teichm\"uller space of a properly embedded
hyperbolic subsurface $\Sigma_S \subseteq \Sigma_{0,5}$ (possibly
with $\Sigma_S = \Sigma_{0,5}$). Let $\Gamma_0$ be the mapping
class group of $\Sigma_S$, and recall that $\Gamma_0$ acts
cocompactly on $\overline{S}$. We will denote $\overline{S}
\setminus S$ by $S_F$. We now have two possibilities:

\noindent {\bf{(a)}} There exists a point $x_0 \in \D$  and a
constant $\delta = \delta(x_0)$ such that $d(\phi_n(x_0), S_F) >
\delta$, for infinitely many $n\in\N$. We can assume, up to
passing to $X/\Gamma_0$ and lifting back, that the sequence
$(\phi_n(x_0))_{n\in\N}$ converges to a point $w \in S$, since the
action of $\Gamma_0$ on $\overline{S}$ is cocompact.

\noindent {\bf {Notation.}} We will write $D_e(a,r)$ to denote the
euclidean disc in $E^2$ with centre $a$ and radius $r$.

\noindent Consider $D_e(x_0,\eta)$ , where $\eta = {\rm{min}}(\delta/4, 1 -
d_e(0,x_0))$. Then,
$$d(\phi_n(x), S_F) \geq d(\phi_n(x_0), S_F) - d(\phi_n(x),
phi_n(x_0)) \geq \delta - 2\eta  \geq \delta/2 >0,$$ for all $x
\in D_e(x_0, \eta)$. So, $d(\phi_n(D_e(x_0,\eta), S_F) >
\delta/2$, for all $n \in \N$. Since $S$ is locally compact, the
maps $\phi_n$ converge (maybe after passing to a subsequence) on
$D_e(x_0,\eta)$ to a map $\phi$. Therefore, $\phi(D_e(x_0,\eta))$
is a copy of a euclidean disc in $S$, which is impossible since
all its sectional curvatures are strictly negative.

\noindent {\bf{(b)}} Otherwise, for all $x \in \D$, $d(\phi_n(x),
S_F) \rightarrow 0$ as $n \rightarrow \infty$. Again we can assume
(up to the action of $\Gamma$) that $\phi_n(0) \rightarrow w'$ as
$n \rightarrow \infty$, for some $w' \in S_F$.

\begin{remark} One of the consequences of the {\it{Collar Lemma}} is
that if $\alpha$ and $\beta$ are intersecting curves on a
hyperbolic surface $\Sigma$, then their lengths cannot be very
small simultaneously. In the light of this result, it is possible
to show (see \cite{Wo}, Corollary $22$) that there exists a
constant $k_0 = k_0(\Sigma)$ such that two strata $S_1,S_2
\subseteq {\overline {\cal{T}}_{WP}(\Sigma) \setminus {\cal{T}}}_
{WP}(\Sigma)$ either have intersecting closures or they satisfy
$d(S_1,S_2) \geq k_0$.
\end{remark}

Let $k_0 = k_0(\Sigma_S)$ be the constant given in the remark
above and consider the maps $\phi_n$ restricted to $D_e(0,k_0/3)$.
We may as well assume that $k_0 \leq 1$ (if not, take $k'_0 =
{\rm{min}} \{1,k_0\}$). Also from the remark above, we deduce that
$S$ must have label $\emptyset$. Otherwise, if $S$ had label
$\alpha$, for some simple closed curve $\alpha$, the fact that
$d(\phi_n(x),S_F) \rightarrow 0$ as $n \rightarrow \infty$, for
all $x \in D_e(0,k_0/3)$ would imply that there exists a curve
$\beta$ in $\Sigma_{0,5}$, disjoint from $\alpha$, such that
$d(\phi_n(x),S_{\alpha \beta}) \rightarrow 0$ as $n \rightarrow
\infty$, for all $x \in D_e(0,k_0/3)$. But we know that a stratum
of the form $S_{\alpha \beta}$ consists of only one point, which
contradicts the construction of the maps $\phi_n$ since we know
that the maps $\phi_n$ contract distances by a factor $\lambda_n$
at most. Thus we assume, from now on, that $S$ has label
$\emptyset$ and so $\overline{S} = X$ and $S_F = X_F$.

We observe that, given $u \in X_F$, there are at most two
non-trivial strata, say $S_{\alpha}$ and $S_{\beta}$, in $X_F$
such that $u \in \overline{S}_{\alpha} \cap \overline{S}_{\beta}$.
These strata correspond to two (possibly equal) simple closed
curves on the surface on which $u$ has nodes. Note that in the
case when there are exactly two strata, we get that $\{u\} =
\overline{S}_{\alpha} \cap \overline{S}_{\beta}$. So we deduce
that $d(\phi_n(x), \overline{S}_{\alpha} \cap
\overline{S}_{\beta}) \rightarrow 0$ as $n \rightarrow \infty$,
for all $x \in D_e(0,k_0/3)$.

Our next aim is to show that there exists $x_1 \in D_e(0, k_0/3)$
and $k_1 = k_1(x_1)>0$, with $k_1 \leq k_0/3$, so that the images
under $\phi_n$ of the points in $D_e(x_1,k_1)$ get uniformly
arbitrarily close to the same stratum as $n \rightarrow \infty$.
Using this result, we will be able to define a distance
non-decreasing projection from $D_e(x_1,k_1)$ to the closure of
that particular stratum.

Suppose that $\phi_n(0) \rightarrow w' \in \overline{S}_{\alpha}
\setminus S_{\alpha \beta}$ as $n \rightarrow \infty$. In
particular, there exists $r>0$ such that, maybe considering a
subsequence, $d(\phi_n(0), S_{\alpha \beta})>r$, say. Let $D =
D(r)$ be the constant given in Corollary $4$ and consider
$D_e(0,k)$, where $k = {\rm{min}} (D/3, k_0/3)$.  We claim that we
can take $(x_1,k_1)$ to be $(0,k)$. Suppose, for contradiction,
that the images of $D_e(0,k)$ do not get uniformly arbitrarily
close to $\overline{S}_{\alpha}$; that is, there exists $K_0>0$, a
subsequence $(\phi_m)_{m\in\N} \subseteq (\phi_n)_{n\in\N}$ and
points $u_m \in \phi_m(D_e(0,k))$ such that $d(u_m,
\overline{S}_{\alpha}) \geq K_0$ for all $m \in \N$. We know that,
for all $m \in \N$, $u_m = \phi_m(x_m)$ for some $x_m \in
D_e(0,k)$ and thus, up to a subsequence, $x_m \rightarrow y \in
D_e(0,k)$. From the construction of the maps $\phi_n$ we have that
$d(u_m,\phi_m(y)) \rightarrow 0$ as $n \rightarrow \infty$.
Therefore, $d(u_m, \overline{S}_{\alpha} \cup
\overline{S}_{\beta}) \rightarrow 0$ as $m \rightarrow \infty$ and
thus $d(u_m, \overline{S}_{\beta}) \rightarrow 0$, since we know
that $d(u_m, \overline{S}_{\alpha}) \geq K_0$ for all $m$. This
represents a contradiction since $d(u_m,\phi_m(0)) \leq D/3$,  the
points $u_m$, and  $\phi_m(0)$ lie in $X \setminus B(r)$, for all
$m \in \N$, but $d(\overline{S}_{\alpha} \setminus (B(r) \cap
\overline{S_\alpha}), \overline{S}_{\beta} \setminus (B(r) \cap
\overline{S_\alpha})) >D$.

The case $\phi_n(0) \rightarrow w' \in S_{\alpha \beta}$ as $n
\rightarrow \infty$ is dealt with in complete analogy, considering
$x_1$ to be any point in $D_e(0,k_0/3) \setminus \{0\}$ and and
defining $k_1$ in a similar way as we did above.

Since $X$ is a $\rm{CAT}(0)$ space and $\overline{S}_{\alpha}$ is
complete and convex we can consider the orthogonal projection
$\pi: X \rightarrow \overline{S}_{\alpha}$ as defined in Theorem
\ref{projection}. Recall that this projection is distance
non-increasing; in particular, for all $x,y \in D_e(x_1,k_1)$ we
have that $$ d(\pi(\phi_n(x)), \pi(\phi_n(y))) \leq d(\phi_n(x),
\phi_n(y)) \leq d_e(x,y).$$

Choose a sequence $(\delta_m)_{m\in\N}$ of positive reals such
that $\delta_m \rightarrow 0$ as $m \rightarrow \infty$. Up to a
subsequence we can assume that, given $m\in\N$, $d(\phi_n(x),
\overline{S}_{\alpha}) <\delta_m$, for $n \geq m$ and for all $x
\in D_e(x_1,k_1)$. We have that $$d(\phi_n(x), \phi_n(y))\leq
d(\phi_n(x), \pi(\phi_n(x))) + d(\phi_n(y), \pi(\phi_n(y)))+
d(\pi(\phi_n(x)), \pi(\phi_n(y))){},$$

\noindent and thus, for all $m\in\N$,

$$d(\pi(\phi_m(x)), \pi(\phi_m(y))) \geq d(\phi_m(x),
\phi_m(y)) -2\delta_m \geq \lambda_m d_e(x,y) -2\delta_m{},$$

Let $\psi_n = \pi\circ\phi_n: D_e(x_1,k_1) \rightarrow
\overline{S}_{\alpha}$. Then exists a sequence
$(\lambda'_n)_{n\in\N}$ of positive real numbers tending to $1$
(we could simply take $\lambda'_n = \lambda_n - 2\delta_n$, since
$k_1 < 1$ and therefore $2\delta_n < 2\delta_n d_e(x,y)$ for all
$x,y \in D_e(x_1,k_1)$) so that
$$\lambda'_n d_e(x,y) \leq d(\psi_n(x), \psi_n(y)) \leq d_e(x,y).$$ We are now back to the situation
described at the beginning of Section $2$, this time in a stratum
of the form $\overline{S}_{\alpha}$. Reasoning in a totally
analogous way to cases {\bf{(a)}} and {\bf{(b)}} we get that
either we can find an isometrically embedded copy of a euclidean
disc in $X$ or else the maps $\phi_n$ collapse  $D_e(x_1,k_1)$ to
a point in $X$, which is impossible since we know that $\phi_n$
decreases distances by a factor $lambda_n$. In any case, we get a
contradiction

\smallskip

\noindent Therefore $X$ is Gromov hyperbolic, as desired.
\end{proof}

\bigskip  
\noindent Javier Aramayona

\noindent Mathematics Institute, University of Warwick, Coventry
CV4 7AL, U.K.

\noindent jaram@maths.warwick.ac.uk


\vspace{1cm}

\begin{thebibliography}{JA04}

\bibitem[Ab]{Ab} W. Abikoff, {\it The real analytic theory of Teichm\"uller space},
Lecture Notes in Mathematics, $820$. Springer, Berlin, $1980$.

\bibitem[Be]{Be} J. Behrstock, {\it{Asymptotic geometry of the mapping class group and Teichm\"uller space}},
arxiv:math.GT/0502367.

\bibitem[Bo]{Bo} B.H. Bowditch, {\it Minkowskian subspaces of non-positively curved metric
spaces} : Bull. London Math. Soc. $27$ ($1995$), no. $6$,
$575$--$584$.

\bibitem[Br]{Br} J. Brock, {\it The Weil-Petersson metric and volumes of 3-dimensional hyperbolic convex cores} :
J. Amer. Math. Soc., $16$ ($2003$), $495$--$535$.

\bibitem[BrFa]{BrFa} J. Brock, B. Farb, {\it Curvature and rank of Teichm\"uller space} : preprint (2001).

\bibitem[Bri]{Bri} M.R. Bridson, {\it{On the existence of flat
planes in spaces of nonpositive curvature}}, Proc. Amer. Math.
Soc. 123 (1995), no. 1, 223--235.

\bibitem[BriHa]{BriHa} M.R. Bridson, A. Haefliger, {\it Metric spaces of non-positive curvature},
Grundlehren der Mathematischen Wissenschaften [Fundamental
Principles of Mathematical Sciences], $319$. Springer-Verlag,
Berlin, $1999$.

\bibitem[DW]{DW} G. Daskalopoulos, R. Wentworth,
{\it{Classification of Weil-Petersson isometries}},  Amer. J.
Math.  125  (2003),  no. 4, 941--975.

\bibitem[Har]{Har} W.J. Harvey, {\it{ Boundary structure of the
modular group}}, Riemann surfaces and related topics: Proceedings
of the 1978 Stony Brook conference. Ann. of Math. Stud. 97,
Princeton, 1981.

\bibitem[Hu]{Hu} Z. Huang, {\it{Asymptotic flatness of the
Weil-Peterssonmetric on Teichm\"uller Space}}, 2001.
arxiv:math.DG/0312419.

\bibitem[Ma]{Ma} H. Masur, {\it{Extension of the
Weil-Petersson metric to the boundary of Teichm\"uller space}},
Duke Math. J., 43(3):623-635, 1976.

\bibitem[Wo]{Wo} S.A. Wolpert, {\it Geometry of the Weil-Petersson completion of Teichm\"uller space},
Surveys in Differential Geometry, VIII: Papers in Honor of Calabi,
Lawson, Siu and Uhlenbeck, editor S. T. Yau, International Press,
Nov. $2003$.

\bibitem[Ya]{Ya} S. Yamada, {\it{Weil-Petersson Completion of
Teichm\"uller Spaces and Mapping Class Group Actions}}, 2001.
arxiv:math.DG/0112001. W. Abikoff, {\it The real analytic theory
of Teichm\"uller space}, Lecture Notes in Mathematics, $820$.
Springer, Berlin, $1980$.

\end{thebibliography}
\end{document}